# RESEARCH   Open Access

# Discrete epidemic models with two time scales

Rafael Bravo de la Parra[1*]  and Luis Sanz-Lorenzo[2]

*Correspondence:
rafael.bravo@uah.es
[1]U.D. Matemáticas, Universidad de Alcalá, Alcalá de Henares, Spain
Full list of author information is available at the end of the article

**Abstract**

The main aim of the work is to present a general class of two time scales discrete-time epidemic models. In the proposed framework the disease dynamics is considered to act on a slower time scale than a second different process that could represent movements between spatial locations, changes of individual activities or behaviors, or others.

To include a sufficiently general disease model, we first build up from first principles a discrete-time susceptible–exposed–infectious–recovered–susceptible (SEIRS) model and characterize the eradication or endemicity of the disease with the help of its basic reproduction number $\mathcal{R}_0$.

Then, we propose a general full model that includes sequentially the two processes at different time scales and proceed to its analysis through a reduced model. The basic reproduction number $\overline{\mathcal{R}}_0$ of the reduced system gives a good approximation of $\mathcal{R}_0$ of the full model since it serves at analyzing its asymptotic behavior.

As an illustration of the proposed general framework, it is shown that there exist conditions under which a locally endemic disease, considering isolated patches in a metapopulation, can be eradicated globally by establishing the appropriate movements between patches.

**Keywords:** Discrete-time epidemic model; Time scales; Disease eradication or persistence

## 1 Introduction

Infectious diseases such as SARS-CoV-2, AIDS, Ebola, or COVID-19 are becoming a part of usual life. They can catastrophically spread and cause a significant number of deaths. It is crucial to use the best management methods to curtail their harmful consequences. The only way to try to compare the effectiveness of these methods is to formulate appropriate mathematical models that help us on making predictions [4]. Mathematical models in epidemiology have a long history of more than two centuries. Most of these models are formulated in continuous time, possibly because of the wealth of analytical tools available for their study. Nevertheless, at least in the last twenty years, time-discrete mathematical epidemics models have been also used with a significant and increasing frequency.

Formulating epidemic models in discrete time has some advantages over the differential equation models, specially when these latter are untractable analytically and one must resort to numerical simulations. Discrete-time models are better implemented in computer simulations when needed, and their parameters can be more easily related to data due to the natural fit of discrete time units to the periodic data collection used in the laboratory or the field [15]. The formulation of discrete-time models should be done directly from first principles and not as a discretization of continuous time models [7, 19]. The mere discretization of continuous models could lead to unfeasible results and usually deviate the focus from the relevant disease dynamics analysis.

Some authors have proposed discrete-time epidemic models without demographic dynamics [1, 6]. The discrete-time epidemic models that we propose in this work follow the literature [13, 27, 29] that does allow for demographic effects. The epidemic and the demographic processes are sequentially included in the model in the following way: at each time interval, first the individuals change their disease status according to the disease flow and then reproduce into the susceptible class and survive following the demographic flow. This constitutes a difference with the case of continuous-time models where these processes are supposed to act instantaneously and simultaneously.

It is frequent in epidemic models that parameters vary by orders of magnitude. In models including disease and demographic dynamics, it is typical that infectious periods have lengths of the order of days, whereas the life spans have them of the order of years [16, 17]. In vector-borne epidemic models in which the vector is an insect, the time scale associated with the vector is usually much faster than the host's one [5, 21, 26]. Individual mobility and behavior can also entail the existence of time scales in the disease dynamics [12, 24]. The models proposed in the previous references are continuous-time and expressed in the form of two time scales models. The analytical tools at disposal for the analysis of this kind of model, quasi-steady-state hypothesis [5, 26] or geometric singular perturbation theory [16, 17, 21, 24, 28], serve at simplifying its analysis by previously reducing the dimension of the system of ordinary differential equations to be studied.

We present in this work a family of time-discrete epidemic models with two time scales. We consider a population divided into groups that we call patches as if they constituted a generalized metapopulation. The patches could represent real spatial locations but also different individual activities or different individual behaviors. The epidemic process acts locally in each patch and may differ from patch to patch. The key point in considering the patches is that slow and fast time scales can be associated with, respectively, the disease dynamics and the process of patch changes that henceforth we call movements.

In [8, 23] it is shown how to construct the kind of discrete-time model with two time scales that we are using. A reduction method that helps to carry out the analytical study of the model is also developed. It is assumed that within a slow time unit the slow process, epidemic, is defined by a map $\mathcal{S}$ and, analogously, within a fast time unit the fast process, movements, is defined by a map $\mathcal{F}$. Thus, the effect of the fast process along a slow time unit can be described by the kth iterate of $\mathcal{F}$, $\mathcal{F}^{(k)}$, where $k$ approximates the time scales ratio. The combined effect of both processes during a slow time unit can then be seen sequentially as the occurrence of $k$ movement episodes followed by a disease dynamics one. In terms of maps $\mathcal{F}$ and $\mathcal{S}$, the associated discrete-time system, denoting by $X$ the population vector and by $t$ the slow time variable, takes the form

$$X(t+1) = \mathcal{S}\bigl(\mathcal{F}^{(k)}\bigl(X(t)\bigr)\bigr). \tag{1}$$

The slow dynamics, represented by $\mathcal{S}$, corresponds to the local disease model in each patch. As it is not possible to use a general disease model, we have proposed as a rather general case a discrete-time SEIRS model based upon the same assumptions of the continuous-time one used in [3]. We have built the model from first principles following the sequential inclusion of the epidemic and the demographic processes as in [13, 27, 29]. The model is the same in all patches, which are distinguished by different parameter values.

In this work the four-dimensional discrete-time system representing the SEIRS model is studied with the help of the basic reproduction number $\mathcal{R}_0$ obtaining sufficient conditions that ensure either the global asymptotic stability of the disease-free equilibrium (DFE) or the uniform persistence of the disease. Thus, disease eradication or endemicity can be straightforwardly characterized.

The fast dynamics, represented by map $\mathcal{F}$, corresponds to the movements of individuals between patches. They are defined for each of the disease compartments ($S$, $E$, $I$, and $R$) by a probability matrix, which in general can depend on the total number of individuals in each compartment across patches.

The proposed full model is represented by the $4n$-dimensional discrete-time system (1), where $n$ is the number of patches. The reduction method [20] applied to the full model leads to a reduced four-dimensional system. If the movement rates are constant, this system is similar to the SEIRS model already studied. The asymptotic analysis of the full model can then be undertaken by studying the asymptotic behavior of this reduced system.

The main aim of this work is to develop a general framework where two time-scale discrete-time epidemic models can be included. This is achieved by proceeding as described above. To complete the presentation, we have also treated a relevant application in a particular setting. We have shown that there exist conditions under which a locally endemic disease, considering every patch in isolation, can be eradicated globally by establishing the appropriate movements.

The structure of the paper is as follows: In Sect. 2, we first propose a discrete-time SEIRS epidemic model, we compute its $\mathcal{R}_0$, and we find conditions for the DFE to be globally asymptotically stable (GAS) when $\mathcal{R}_0 < 1$ and for the disease to persist when $\mathcal{R}_0 > 1$. Section 3 is devoted to the presentation and reduction of a model combining the disease dynamics at a slow time scale with a fast process interpretable as movements between patches. The previous general model is analyzed in Sect. 4 in the particular case of constant recruitment, standard incidence, and constant movement rates in order to illustrate the advantages of the reduction procedure. We summarize our results and perspectives in the final Discussion section.

## 2 Discrete SEIRS model

In this section we present a discrete-time SEIRS epidemic model which is a discrete version of the continuous SEIRS model included in [3].

The individuals of the population are divided into four epidemiological classes that, at each time $t \in \{0, 1, 2, \dots\}$, are represented by $S(t)$, susceptible, $E(t)$, exposed, $I(t)$, infectious, and $R(t)$, recovered. So, the state vector of the population is $X(t) = (S(t), E(t), I(t), R(t))$, and the total population $N(t) = S(t) + E(t) + I(t) + R(t)$.

Following [27, 29], we consider that in each time interval there exist two distinct temporal phases. In the first one, the disease dynamics acts and individuals can change from one epidemiological class to the next one, as shown in Fig. 1. Reproduction and survival happen in the second temporal phase. We assume that the disease does not affect the birth process. Concerning survival, to take into account both the natural and the disease induced mortalities in full generality, we associate different survival rates with each of the epidemiological classes. The corresponding diagram is included in Fig. 2.

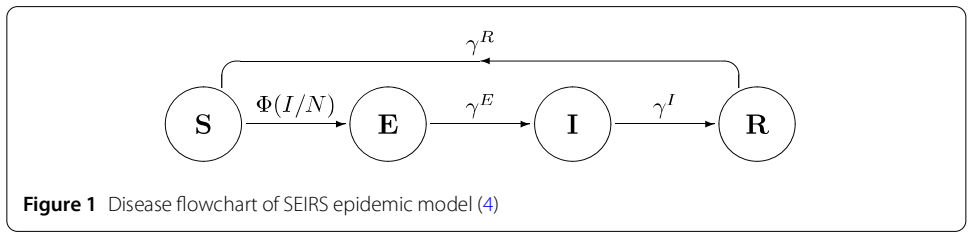

**Figure 1** Disease flowchart of SEIRS epidemic model (4)

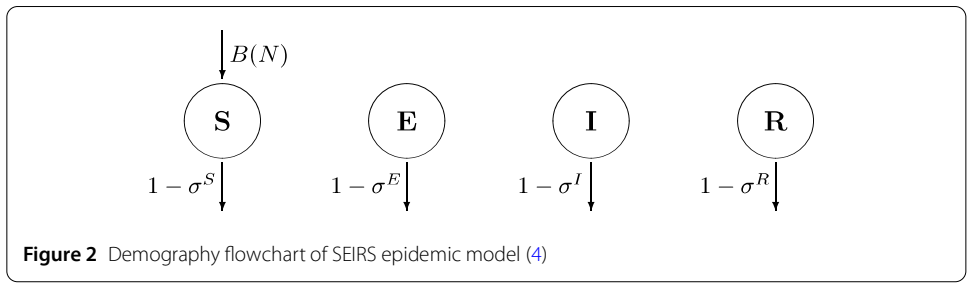

**Figure 2** Demography flowchart of SEIRS epidemic model (4)

The disease transmission is a function $\Phi(I(t)/N(t))$ of the fraction of susceptible individuals that become exposed, where

$$\Phi : [0,1] \longrightarrow [0,1], \qquad \Phi \in C^2([0,1]), \qquad \Phi(0) = 0, \quad \text{and is increasing.} \tag{2}$$

A common choice, which corresponds to the so-called proportional or standard incidence, is $\Phi(x) = \beta x$, where $\beta \in (0,1]$ is the transmission parameter. When infections are modeled as Poisson processes, then $\Phi(x) = 1 - e^{-\beta x}$ for $\beta > 0$ [29].

Upon infection, the transitions between classes are defined by parameters $\gamma^C \in (0,1)$ for $C \in \{E, I, R\}$ that represent the fraction of the individuals of class $C$ that pass to the next class per time unit (see Fig. 1).

Henceforth, we use $C$ as a generic letter for an unspecified epidemiological class and call $\mathcal{C} := \{S, E, I, R\}$.

Demography is included in the model in a simple form. It is assumed that there is no vertical transmission of the disease, so that all births occur into the susceptible class. The recruitment of individuals to the susceptible class per time unit is a function

$$B : [0, \infty) \longrightarrow [0, \infty), \quad B \in C^1([0, \infty)), \tag{3}$$

of the total population $N$. We will use a constant recruitment function as a particular case. Other common choices, proposed in [27, 29], are geometric, Beverton–Holt, and Ricker recruitment functions.

Individuals in all epidemiological classes are subject to natural death, possibly affected by the disease. Parameters $\sigma^C \in (0,1)$, for $C \in \mathcal{C}$, represent the fraction of the individuals of class $C$ that survive per time unit. Thus, the fraction of individuals dying per time unit is $1 - \sigma^C \in (0,1)$ (see Fig. 2).

The transitions between epidemiological classes in a time unit are defined by the following map:

$$T(X) = \begin{pmatrix} S - \Phi(I/N)S + \gamma^R R \\ E + \Phi(I/N)S - \gamma^E E \\ I + \gamma^E E - \gamma^I I \\ R + \gamma^I I - \gamma^R R \end{pmatrix}$$

and the demographic changes by the next one

$$D(X) = \left(B(N) + \sigma^S S, \sigma^E E, \sigma^I I, \sigma^R R\right).$$

We now take into account both processes sequentially: first epidemiological transitions followed by demography. The discrete-time SEIRS epidemic model that we propose can so be expressed in terms of maps $T$ and $D$ as

$$X(t+1) = D\bigl(T\bigl(X(t)\bigr)\bigr),$$

or in a detailed form

$$\begin{aligned}
S(t+1) &= B\bigl(N(t)\bigr) + \sigma^S \gamma^R R(t) + \sigma^S \left(1 - \Phi\left(\frac{I(t)}{N(t)}\right)\right) S(t), \\
E(t+1) &= \sigma^E \Phi\left(\frac{I(t)}{N(t)}\right) S(t) + \sigma^E \bigl(1 - \gamma^E\bigr) E(t), \\
I(t+1) &= \sigma^I \gamma^E E(t) + \sigma^I \bigl(1 - \gamma^I\bigr) I(t), \\
R(t+1) &= \sigma^R \gamma^I I(t) + \sigma^R \bigl(1 - \gamma^R\bigr) R(t).
\end{aligned} \quad (4)$$

The assumptions on the parameters of the model allow us to straightforwardly prove that $T(\mathbf{R}_+^4) \subset \mathbf{R}_+^4$ and also $D(\mathbf{R}_+^4) \subset \mathbf{R}_+^4$. Therefore, $\mathbf{R}_+^4$ is forward invariant under the semi-flow defined by the discrete-time system (4).

**Proposition 1** *If the function B is bounded, then system* (4) *is dissipative.*

*Proof of Proposition* 1 Let $\hat{B} > 0$ be such that $B(x) \leq \hat{B}$ for all $x \in [0, \infty)$, and $\hat{\sigma} = \max_{C \in \mathcal{C}}\{\sigma^C\} \in (0, 1)$. Then

$$\begin{aligned}
N(t+1) &= B\bigl(N(t)\bigr) + \left(\sigma^S\left(1 - \Phi\left(\frac{I(t)}{N(t)}\right)\right) + \sigma^E \Phi\left(\frac{I(t)}{N(t)}\right)\right) S(t) \\
&\quad + \bigl(\sigma^E(1-\gamma^E) + \sigma^I \gamma^E\bigr) E(t) + \bigl(\sigma^I(1-\gamma^I) + \sigma^R \gamma^I\bigr) I(t) \\
&\quad + \bigl(\sigma^R(1-\gamma^R) + \sigma^S \gamma^R\bigr) R(t) \\
&\leq \hat{B} + \max\{\sigma^S, \sigma^E\} S(t) + \max\{\sigma^E, \sigma^I\} E(t) \\
&\quad + \max\{\sigma^I, \sigma^R\} I(t) + \max\{\sigma^R, \sigma^S\} R(t) \\
&\leq \hat{\sigma} N(t) + \hat{B}.
\end{aligned} \quad (5)$$

We now present some basic properties for the solution of the linear scalar difference equation

$$x(t+1) = ax(t) + b \qquad (6)$$

with $0 < a < 1$, $b \geq 0$ that will be used here and later on in the manuscript. The solution to (6) is

$$x(t) = \left(x(0) - \frac{b}{1-a}\right)a^t + \frac{b}{1-a}, \qquad (7)$$

from where it follows that $x(t)$ converges monotonically to $b/(1-a)$. In particular, for all $x(0) \geq 0$, we have

$$\min\left\{x(0), \frac{b}{1-a}\right\} \leq x(t) \leq \max\left\{x(0), \frac{b}{1-a}\right\}, \quad t = 0, 1, 2, \ldots \qquad (8)$$

Therefore, from (5) and (7) it follows that any solution of system (4), with initial conditions $X(0) \in \mathbf{R}_+^4$, satisfies that

$$N(t) \leq \left(N(0) - \frac{\hat{B}}{1-\hat{\sigma}}\right)\hat{\sigma}^t + \frac{\hat{B}}{1-\hat{\sigma}} \xrightarrow[t\to\infty]{} \frac{\hat{B}}{1-\hat{\sigma}}.$$

This inequality proves that, fixed any $M > \hat{B}/(1-\hat{\sigma})$, for every solution of (4), there is a time $\hat{t}$ such that $N(t) \leq M$ for all $t \geq \hat{t}$ and, thus, the dissipativity of system (4). □

Function $B$ is bounded in the particular case of having a constant, Beverton–Holt, or Ricker recruitment function.

Notice that in the proof it is shown that all nonnegative solutions of system (4) are attracted by the compact set $K = \{(S, E, I, R) \in \mathbf{R}_+^4 : N \in [0, \hat{B}/(1-\hat{\sigma})]\}$.

We continue the analysis by finding conditions for the eradication or endemicity of the disease.

To consider disease eradication, we try to find an equilibrium of (4) with $I = 0$. We immediately obtain $E = 0$, $R = 0$, and $S = B(S) + \sigma^S S$. To guarantee the existence of a unique disease-free equilibrium (DFE), we make the following assumption on the scalar difference equation representing the demography of the population without disease.

**Hypothesis 2.1** *The equation $S(t+1) = \sigma^S S(t) + B(S(t))$ possesses a unique positive equilibrium $S^*$ that is hyperbolic and globally asymptotically stable (GAS) in $(0, \infty)$.*

If the recruitment function is constant, then Hypothesis 2.1 holds. This also happens, for certain values of the parameters, in the cases of the Beverton–Holt and the Ricker recruitment functions.

Note that if Hypothesis 2.1 is met, then the unique DFE of system (4) is $X_0^* = (S^*, 0, 0, 0)$.

As we will prove, the basic reproduction number $\mathcal{R}_0$ of system (4) determines whether the disease is eradicated or it becomes endemic. We use the next-generation method to calculate $\mathcal{R}_0$ of a discrete-time system as developed in [2]. We consider as infected states $E$ and $I$, and as uninfected states $S$ and $R$. The so-called disease-free system is the one associated with the uninfected states setting $E = 0$ and $I = 0$:

$$S(t+1) = B(S(t) + R(t)) + \sigma^S \gamma^R R(t) + \sigma^S S(t),$$
$$R(t+1) = \sigma^R (1 - \gamma^R) R(t). \tag{9}$$

**Proposition 2** *If function $B$ is bounded and Hypothesis* 2.1 *is met, then $(S^*, 0)$ is a GAS equilibrium of system* (9) *in $\Omega = \{(S, R) \in \mathbf{R}_+^2 : S > 0\}$.*

*Proof of Proposition* 2  It is straightforward to show that $(S^*, 0)$ is an equilibrium of system (9) and, by linearization, that it is hyperbolic and locally asymptotically stable (LAS).

To prove that $(S^*, 0)$ attracts all points in $\Omega$, we apply Theorem 2.1 in [30].

For that, let $(S_0, R_0) \in \Omega$ and define, for $t = 0, 1, 2, \ldots$,

$$\sigma_t(S) = B\big(S + \big(\sigma^R(1 - \gamma^R)\big)^t R_0\big) + \sigma^S \gamma^R \big(\sigma^R(1 - \gamma^R)\big)^t R_0 + \sigma^S S,$$

that is, a continuous map in $\mathbf{R}_+$, and the discrete dynamical process defined by $\tau_0 := I$, the identity map, and $\tau_t := \sigma_{t-1} \circ \sigma_{t-2} \circ \cdots \circ \sigma_1 \circ \sigma_0$ for $t \geq 1$.

If $\{(S(t), R(t)) : t \geq 0\}$ is the orbit associated with $(S_0, R_0)$ in system (9), we have that $R(t) = (\sigma^R(1 - \gamma^R))^t R_0$ and $S(t) = \tau_t(S_0)$.

Since $0 < \sigma^R(1 - \gamma^R) < 1$, we have $\lim_{t \to \infty} R(t) = 0$, and so we need to show that $\lim_{t \to \infty} S(t) = S^*$.

We now prove that the discrete dynamical process $\tau_t$ ($t \geq 0$) is asymptotically autonomous (Definition 2.1 in [30]) and that its limit discrete semiflow, $\Sigma_t = \Sigma \circ \overset{(t)}{\cdots} \circ \Sigma$ ($t \geq 0$), is that generated by the continuous map $\Sigma(S) = B(S) + \sigma^S S$ in $\mathbf{R}_+$. In order to do so, let $x \in \mathbf{R}_+$, and let $\{x_n\}_{n=1}^\infty$ and $\{t_n\}_{n=1}^\infty$ be sequences in $\mathbf{R}_+$ such that $\lim_{n \to \infty} x_n = x$ and $\lim_{n \to \infty} t_n = \infty$. We need to show that $\lim_{n \to \infty} \Sigma_{t_n}(x_n) = \Sigma(x)$, which follows immediately taking into account that $0 < \sigma^R(1 - \gamma^R) < 1$ and that $B$ is continuous.

As the only $\Sigma$-invariant subset of $\mathbf{R}_+$ is $\{S^*\}$, to complete the proof by applying Theorem 2.1 in [30], we just need to prove that the set $\{S(t) : t \geq 0\}$ is bounded.

Let $\hat{B}$ be un upper-bound of function $B$. Then, for every $t \geq 0$, we have $\sigma_t(S) \leq \hat{B} + R_0 + \sigma^S S$, and so $\sigma_t(S_0) \leq x(t)$, where $x(t)$ is the solution to

$$x(t+1) = \hat{B} + R_0 + \sigma^S x(t),$$
$$x(0) = S_0,$$

which corresponds to the linear scalar equation (6) with $a := \sigma^S$ and $b := \hat{B} + R_0$. Now, using (8), we have $\sigma_t(S_0) \leq x(t) \leq \max\{S_0, (\hat{B} + R_0)/(1 - \sigma^S)\}$, and so

$$\sup\{S(t) : t \geq 0\} = \sup\{\tau_t(S_0) : t \geq 0\} \leq \max\{S_0, (\hat{B} + R_0)/(1 - \sigma^S)\},$$

as we wanted to show. □

To proceed with the application of the next-generation method, in the equations for the infected compartments we must separate the terms $\mathcal{F}_E$ and $\mathcal{F}_I$, representing new infec-tions, from the terms $\mathcal{T}_E$ and $\mathcal{T}_I$ associated with transitions between compartments:

$$\mathcal{F}_E(X) = \sigma^E \Phi(I/N)S, \qquad \mathcal{F}_I(X) = 0,$$
$$\mathcal{T}_E(X) = \sigma^E(1-\gamma^E)E, \qquad \mathcal{T}_I(X) = \sigma^I \gamma^E E + \sigma^I(1-\gamma^I)I.$$

Now, we can calculate matrices

$$F = \left[\frac{\partial \mathcal{F}_C(X_0^*)}{\partial D}\right]_{C,D \in \{E,I\}} = \begin{bmatrix} 0 & \sigma^E \Phi'(0) \\ 0 & 0 \end{bmatrix}$$

and

$$T = \left[\frac{\partial \mathcal{T}_C(X_0^*)}{\partial D}\right]_{C,D \in \{E,I\}} = \begin{bmatrix} \sigma^E(1-\gamma^E) & 0 \\ \sigma^I \gamma^E & \sigma^I(1-\gamma^I) \end{bmatrix}$$

to obtain the next-generation matrix $Q = F(Id - T)^{-1}$ [2], where $Id$ is the identity matrix of order 2 whose spectral radius is $\mathcal{R}_0$.

$$\mathcal{R}_0 = \rho(Q) = \frac{\sigma^E \sigma^I \gamma^E \Phi'(0)}{(1 - \sigma^E(1-\gamma^E))(1 - \sigma^I(1-\gamma^I))}. \tag{10}$$

The next result gives sufficient conditions for the disease eradication locally and globally.

**Theorem 3** *Let system* (4) *satisfy Hypothesis* 2.1. *Then*
  (a) *If* $\mathcal{R}_0 < 1$, *then DFE* $X_0^*$ *is LAS.*
  (b) *If* $\mathcal{R}_0 < 1$ *and* $\Phi''(x) \leq 0$ *for* $x > 0$, *then DFE* $X_0^*$ *is GAS in*
      $\Omega = \{(S, E, I, R) \in \mathbf{R}_+^4 : S > 0\}$.
  (c) *If* $\mathcal{R}_0 > 1$, *then DFE* $X_0^*$ *is unstable.*

*Proof of Theorem* 3
  (a) and (c) are direct consequences of Theorem 2.1 in [2].
  (b) $\Phi \in C^2([0,1])$, $\Phi(0) = 0$, and $\Phi''(x) \leq 0$ imply that $\Phi(x) \leq \Phi'(0)x$ for $x \in [0,1]$.
Therefore

$$E(t+1) \leq \sigma^E \Phi'(0)\frac{I(t)}{N(t)}S(t) + \sigma^E(1-\gamma^E)E(t)$$
$$\leq \sigma^E \Phi'(0)I(t) + \sigma^E(1-\gamma^E)E(t),$$

and so, using the equation for $I$ in (4),

$$\begin{bmatrix} E(t+1) \\ I(t+1) \end{bmatrix} \leq \begin{bmatrix} \sigma^E(1-\gamma^E) & \sigma^E \Phi'(0) \\ \sigma^I \gamma^E & \sigma^I(1-\gamma^I) \end{bmatrix} \begin{bmatrix} E(t) \\ I(t) \end{bmatrix}$$
$$= (F+T)\begin{bmatrix} E(t) \\ I(t) \end{bmatrix} \leq (F+T)^{t+1} \begin{bmatrix} E(0) \\ I(0) \end{bmatrix}.$$

Matrix $F + T$ can be considered to be the projection matrix of a standard linear matrix model of population dynamics [18], with matrices $F$ and $T$ representing the fertility and transition matrices respectively. The net reproductive rate of the model coincides with $\mathcal{R}_0$ and, since $\rho(T) = \max\{\sigma^E(1-\gamma^E), \sigma^I(1-\gamma^I)\} < 1$ and $\mathcal{R}_0 < 1$, we can apply Theorem 3.3 in [18], which yields that $\rho(F+T) < 1$ and, therefore, $(F+T)^t \underset{t\to\infty}{\longrightarrow} 0$. This proves that $E(t)$, $I(t) \underset{t\to\infty}{\longrightarrow} 0$.

To prove that also $R(t)$ tends to 0, we use from the previous arguments that we can find $\alpha \in (0,1)$ and $K > 0$ such that $I(t) \leq K\alpha^t$. Thus, substituting in the equation of $R(t)$, we obtain $R(t+1) \leq K\alpha^t + \sigma^R(1-\gamma^R)R(t)$, which implies that there exist $\bar{\alpha} \in (\max(\alpha, \sigma^R(1-\gamma^R)), 1)$ and $\bar{K} > 0$ such that $R(t) \leq \bar{K}\bar{\alpha}^t \underset{t\to\infty}{\longrightarrow} 0$.

Finally, to prove that $S(t) \underset{t\to\infty}{\longrightarrow} S^*$, we just need to follow the reasoning in the proof of Proposition 2. Defining, for $(S_0, E_0, I_0, R_0) \in \Omega$ and $t = 0, 1, 2, \ldots$,

$$\sigma_t(S) = B\bigl(S + E(t) + I(t) + R(t)\bigr) + \sigma^S \gamma^R R(t)$$
$$+ \sigma^S\left(1 - \Phi\left(\frac{I(t)}{S + E(t) + I(t) + R(t)}\right)\right)S.$$

The corresponding discrete dynamical process $\tau_t$ ($t \geq 0$) is asymptotically autonomous with limit discrete semiflow $\Sigma^t$ ($t \geq 0$) generated by the continuous map $\Sigma(S) = B(S) + \sigma^S S$ in $\mathbf{R}_+$, which coincides with the one in the proof of Proposition 2. □

Note that the assumption $\Phi''(x) \leq 0$ for $x > 0$ is met in the case of standard incidence, i.e., if $\Phi(x) = \beta x$.

The endemicity of the disease is represented in mathematical terms by the concept of uniform persistence. We use the persistence function $\rho(S, E, I, R) = E + I$. Thus, system (4) is *uniformly persistent* [25] if there exists $\varepsilon > 0$ such that $\liminf_{t\to\infty}(E(t) + I(t)) > \varepsilon$ for any solution with $E(0) + I(0) > 0$. If lim inf is substituted by lim sup in the definition, the system is said to be *uniformly weakly persistent*.

In the next theorem we prove the uniform persistence of system (4) when $\mathcal{R}_0 > 1$ in the case of constant recruitment function $B(N) = B$ and standard incidence $\Phi(x) = \beta x$, $\beta \in [0,1)$.

**Theorem 4** *Let $\Phi(x) = \beta x$, $0 < \beta \leq 1$, and $B(N) = B$ be constant in system (4). If $\mathcal{R}_0 > 1$, then (4) is uniformly persistent.*

*Proof of Theorem* 4 Since all nonnegative solutions of (4) are attracted by a compact set and $\{X \in \mathbf{R}_+^4 : E + I > 0\}$ is forward invariant, Corollary 4.8 in [25] establishes that it is sufficient to prove that it is uniformly weakly persistent to obtain that it is also uniformly persistent.

So, let us prove that (4) is uniformly weakly persistent. We argue by contradiction. Suppose that it is not. Then, for any arbitrary $\varepsilon > 0$, there exists a solution $X(t)$ with $E(0) + I(0) > 0$ and $\limsup_{t\to\infty}(E(t) + I(t)) < \varepsilon$. Thus, there exists some $t_0 > 0$ such that $E(t) + I(t) < \varepsilon$ for all $t \geq t_0$.

Also, for $t \geq t_0$, $R(t+1) < \sigma^R(1-\gamma^R)R(t) + \sigma^R\gamma^I\varepsilon$ with $0 < \sigma^R(1-\gamma^R) < 1$. This implies, iterating the right-hand side from $t_0$ on,

$$R(t) < \left(\sigma^R(1-\gamma^R)\right)^{t-t_0} R(t_0) + \sigma^R\gamma^I\varepsilon \sum_{i=0}^{t-t_0-1} \left(\sigma^R(1-\gamma^R)\right)^i$$

$$< \left(\sigma^R(1-\gamma^R)\right)^{t-t_0} R(t_0) + \frac{\sigma^R\gamma^I}{(1-\sigma^R)+\gamma^R-(1-\sigma^R)\gamma^R}\varepsilon.$$

Then, as $(\sigma^R(1-\gamma^R))^{t-t_0} \xrightarrow[t\to\infty]{} 0$, there exists $t_1 > t_0$ such that $R(t) < (1 + \frac{\sigma^R\gamma^I}{\gamma^R+d^R-\gamma^R d^R})\varepsilon := \bar{r}\varepsilon$ for $t \geq t_1$.

Thus, for any $\varepsilon > 0$, there exists $t_1 > 0$ such that

$$E(t) + I(t) + R(t) \leq (1+\bar{r})\varepsilon \quad \text{for } t \geq t_1. \tag{11}$$

Let us now establish a lower bound for the total population $N(t)$. Reproducing, with the appropriate changes, the calculations in the proof of Proposition 1, we obtain the following inequality:

$$N(t+1) \geq \bar{\sigma}N(t) + B,$$

where $\bar{\sigma} = \min_{C\in\mathcal{C}}\{\sigma^C\} \in (0,1)$ that implies, considering the difference equation $x(t+1) = \bar{\sigma}x(t) + B$ and using (8),

$$N(t) \geq \bar{n} := \min\{N(0), B/(1-\bar{\sigma})\} > 0.$$

This inequality together with (11) yields

$$S(t) \geq \bar{n} - (1+\bar{r})\varepsilon \quad \text{for } t \geq t_1. \tag{12}$$

With the help of (11) and (12), we can find the following lower bound for the coefficient of $I(t)$ in the $E$ equation:

$$\frac{\sigma^E\beta S(t)}{N(t)} \geq \frac{\sigma^E\beta S(t)}{S(t)+(1+\bar{r})\varepsilon} \geq \frac{\sigma^E\beta(\bar{n}-(1+\bar{r})\varepsilon)}{\bar{n}-(1+\bar{r})\varepsilon+(1+\bar{r})\varepsilon}$$

$$= \sigma^E\beta\left(1 - \frac{1+\bar{r}}{\bar{n}}\varepsilon\right) =: G(\varepsilon).$$

Let us define matrix

$$\bar{P}_\varepsilon = \begin{bmatrix} \sigma^E(1-\gamma^E) & G(\varepsilon) \\ \sigma^I\gamma^E & \sigma^I(1-\gamma^I) \end{bmatrix}$$

that satisfies, for $t \geq t_1$,

$$\begin{bmatrix} E(t+1) \\ I(t+1) \end{bmatrix} \geq \bar{P}_\varepsilon \begin{bmatrix} E(t) \\ I(t) \end{bmatrix}.$$

As $\bar{P}_\varepsilon$ is a primitive matrix, if we can find $\varepsilon$ such that $\rho(\bar{P}_\varepsilon) > 1$, we obtain that $E(t), I(t) \underset{t\to\infty}{\longrightarrow} \infty$ whenever $E(0) + I(0) > 0$, which is the contradiction we were looking for. As seen in the proof of Theorem 3, we can check it by means of the net reproductive rate $\mathcal{R}_{0,\varepsilon}$ of the model associated with matrix $\bar{P}_\varepsilon$ (Theorem 3.3. in [18]).

Now $\mathcal{R}_{0,\varepsilon}$ and $\mathcal{R}_0$ satisfy

$$\mathcal{R}_{0,\varepsilon} = \frac{G(\varepsilon)}{\sigma^E \beta}\mathcal{R}_0 = \left(1 - \frac{1+\bar{r}}{\bar{n}}\varepsilon\right)\mathcal{R}_0.$$

So, if we choose $\bar{\varepsilon} = \frac{1}{2}(1 - \frac{1}{\mathcal{R}_0})\frac{\bar{n}}{1+\bar{r}}$, and take ito account that $\mathcal{R}_0 > 1$, we obtain the required result to complete the proof

$$\mathcal{R}_{0,\bar{\varepsilon}} = \frac{1}{2}(\mathcal{R}_0 + 1) > 1. \qquad \square$$

## 3 The model

In this section we present a model of disease dynamics with two time scales. The population can be considered divided into groups that we call patches as if they were forming a sort of *generalized metapopulation*. The patches could represent real spatial locations and involve explicit movements of the individuals between them [3], but also different individual daily activities (places of residency, work, or business) [14] or different individual behaviors [12]. The key point in considering the *generalized patches* is that the disease dynamics occurs locally and at a slow time scale compared to the fast time scale associated with the *movements* between patches.

We consider that individuals move between $n$ patches. In each patch the disease dynamics follows SEIRS model (4) with appropriate local parameters. We assume that movements are almost instantaneous with respect to the disease dynamics. Thus, the model takes the form of a time-discrete two time scale system (1) with movements being the fast process and the disease dynamics the slow process [8].

We denote the densities of susceptible, exposed, infectious, and recovered individuals in patch $j \in \{1,\ldots,n\}$ at time $t \in \{0,1,2,\ldots\}$ by $S_j(t)$, $E_j(t)$, $I_j(t)$, and $R_j(t)$, respectively, and the total population $N_j(t) = S_j(t) + E_j(t) + I_j(t) + R_j(t)$.

We denote by $\bar{x}_C = \mathrm{col}(C_1,\ldots,C_n) \in \mathbb{R}^{n\times 1}$ for $C \in \mathcal{C}$ ($\mathcal{C} = \{S,E,I,R\}$) the state vectors of individuals in each compartment (susceptible, exposed, infective, and recovered) across the $n$ patches. The population, or rather metapopulation, state vector is called

$$X = \mathrm{col}(\bar{x}_S, \bar{x}_E, \bar{x}_I, \bar{x}_R) \in \mathbb{R}^{4n\times 1}.$$

The existence of two time scales in the complete model that we propose leads to a reduced model for some global variables. In this case the global variables correspond to the total number of individual in each compartment:

$$S = \sum_{j=1}^n S_j, \qquad E = \sum_{j=1}^n E_j, \qquad I = \sum_{j=1}^n I_j, \qquad R = \sum_{j=1}^n R_j,$$

that we collect in the vector of global variables

$$Y = \mathrm{col}(S,E,I,R),$$

whose sum yields the total number of individuals in the metapopulation

$$N = S + E + I + R = \sum_{j=1}^{n} N_j.$$

It is straightforward to see that we can obtain the global variables from the state variables with the help of matrix $U = \text{diag}(\bar{1}, \bar{1}, \bar{1}, \bar{1}) \in \mathbb{R}_+^{4 \times 4n}$, where $\bar{1} = (1, \overset{(n)}{\ldots}, 1) \in \mathbb{R}_+^n$ is a row vector

$$Y = UX.$$

*Fast process: movements*   We assume that individuals in each compartment move between patches according to movement rates that can generally depend on the global variables $Y$. In this way, for each disease compartment, movements are represented by a regular stochastic matrix depending on $Y \in \mathbb{R}_+^4$

$$M^S(Y), M^E(Y), M^I(Y), M^R(Y) \in \mathbb{R}_+^{n \times n}.$$

The movements of the whole metapopulation are then defined through the following matrix:

$$M(Y) = \text{diag}\left(M^S(Y), M^E(Y), M^I(Y), M^R(Y)\right) \in \mathbb{R}_+^{4n \times 4n}.$$

The state $X$ of the metapopulation after one movement episode is defined by the following map:

$$\mathcal{F}(X) = M(UX)X \tag{13}$$

that represents the fast process in system (1).

*Slow process: disease dynamics*   The slow process, the disease dynamics, is defined locally, i.e., in each patch $j \in \{1, \ldots, n\}$, by SEIRS model (4):

$$\begin{aligned}
S_j(t+1) &= B_j\left(N_j(t)\right) + \sigma_j^S \gamma_j^R R_j(t) + \sigma_j^S \left(1 - \Phi_j\left(\frac{I_j(t)}{N_j(t)}\right)\right) S_j(t), \\
E_j(t+1) &= \sigma_j^E \Phi_j\left(\frac{I_j(t)}{N_j(t)}\right) S_j(t) + \sigma_j^E \left(1 - \gamma_j^E\right) E_j(t), \\
I_j(t+1) &= \sigma_j^I \gamma_j^E E_j(t) + \sigma_j^I \left(1 - \gamma_j^I\right) I_j(t), \\
R_j(t+1) &= \sigma_j^R \gamma_j^I I_j(t) + \sigma_j^R \left(1 - \gamma_j^R\right) R_j(t).
\end{aligned} \tag{14}$$

The recruitment functions $B_j$ verify assumption (3), and the transmission functions $\Phi_j$ assumption (2). All the parameters $\sigma_j^C$ and $\gamma_j^C$ are in $(0,1)$, $j \in \{1, \ldots, n\}$ and $C \in \mathcal{C}$.

To obtain the map $\mathcal{S}$ representing the slow process, we need to appropriately reorder equations (14) for all $n$ patches. Following the order of variables in the population vector $X$, we must include first the equations for variables $S_1, \ldots, S_n$ and then follow consecutively with those corresponding to compartments $E$, $I$, and $R$.

Finally, the complete two time scale model takes the form of system (1)

$$X(t+1) = \mathcal{S}\left(\mathcal{F}^{(k)}(X(t))\right).$$

Note that

$$UF(X) = UM(UX)X = UX,$$

and so map $\mathcal{F}$ keeps invariant the values of the global variables $Y$ and, therefore, its kth iterate can be expressed in terms of the k-power $M(Y)^k$ of matrix $M(Y)$

$$\mathcal{F}^{(k)}(X) = M(UX)^k X,$$

so that the complete model reads as follows:

$$X(t+1) = \mathcal{S}\left(M(UX(t))^k X(t)\right). \tag{15}$$

**Proposition 5** *If functions $B_j$ are bounded, $j \in \{1, \ldots, n\}$, then system* (15) *is dissipative.*

*Proof of Proposition* 5 Let us represent with hat $\hat{C}_j(t+1)$ ($C \in \mathcal{C}$, $j \in \{1, \ldots, n\}$) the state variables after the fast process and before the disease dynamics acts.

We know that $\sum_{j=1}^n \hat{C}_j(t+1) = \sum_{j=1}^n C_j(t)$ for $C \in \mathcal{C}$.

Moreover, following the proof of Proposition 1, we have, for every $j \in \{1, \ldots, n\}$, that

$$\sum_{C \in \mathcal{C}} C_j(t+1) \leq \bar{\sigma}_j \sum_{C \in \mathcal{C}} \hat{C}_j(t+1) + \hat{B}_j,$$

where $\hat{B}_j > 0$ is such that $B_j(x) \leq \hat{B}_j$ for all $x \in [0, \infty)$, and $\hat{\sigma}_j = \max_{C \in \mathcal{C}} \{\sigma_j^C\} \in (0, 1)$. Now, calling $\hat{B} = \hat{B}_1 + \cdots + \hat{B}_n$ and $\hat{\sigma} = \max_{j \in \{1, \ldots, n\}} \{\hat{\sigma}_j\}$, we obtain the following recurrent inequality for the total population:

$$\begin{aligned} N(t+1) &= \sum_{j=1}^n \left(\sum_{C \in \mathcal{C}} C_j(t+1)\right) \\ &\leq \sum_{j=1}^n \left(\hat{\sigma}_j \sum_{C \in \mathcal{C}} \hat{C}_j(t+1)\right) + \sum_{j=1}^n \hat{B}_j \\ &\leq \hat{\sigma} \sum_{j=1}^n \left(\sum_{C \in \mathcal{C}} \hat{C}_j(t+1)\right) + \hat{B} = \hat{\sigma} N(t) + \hat{B}. \end{aligned}$$

Thus, the proof can be completed as in Proposition 1. □

Notice that all nonnegative solutions of system (15) are attracted by the compact set $K = \{X \in \mathbf{R}_+^{4n} : N \in [0, \bar{B}/(1 - \bar{\sigma})]\}$.

### 3.1 Reduced system

The fact that the probability matrix $M^C(Y)$ is primitive for every $C \in \mathcal{C}$ and $Y \in \mathbb{R}_+^4$ implies that 1 is its strictly dominant eigenvalue, $\bar{1}$ is an associated row left eigenvector, and there exists a unique column right eigenvector $\bar{m}^C(Y)$, representing the corresponding stable probability distribution, that satisfies $\bar{1}\bar{m}^C(Y) = 1$.

In the reduction procedure of model (15) we need to calculate the limit of the iterates of map $\mathcal{F}$ that, in this case, is equivalent to calculating the limit of the powers of matrix $M(Y)$. This latter follows from the Perron–Frobenius theorem:

$$\lim_{k\to\infty} M^C(Y)^k = \bar{m}^C(Y)\bar{1} \quad \text{and} \quad \lim_{k\to\infty} M(Y)^k = \bar{M}(Y)U,$$

where $\bar{M}(Y) = \text{diag}(\bar{m}^S(Y), \bar{m}^E(Y), \bar{m}^I(Y), \bar{m}^R(Y))$. Finally,

$$\lim_{k\to\infty} \mathcal{F}^{(k)}(X) = \bar{\mathcal{F}}(X) := \bar{M}(UX)UX. \tag{16}$$

Thus, the four-dimensional reduced model to be used to study the asymptotic behavior of the solutions of system (15) is

$$Y(t+1) = U\mathcal{S}(\bar{M}(Y(t))Y(t)). \tag{17}$$

The next result states how the analysis of the stability of the equilibria of system (17) extends to system (15). It is a direct translation to our setting of Theorem 2 in [20]. Loosely speaking, for $k$ large enough, from the equilibria of the reduced system we can obtain good approximations of the equilibria of system (15) and, in case of asymptotic stability, of their basins of attraction.

**Theorem 6** *Let $Y^* \in \mathbb{R}_+^4$ be a hyperbolic equilibrium point of system* (17). *Then there exists an integer $k_0 \geq 0$ such that, for all $k \geq k_0$, system* (15) *has an equilibrium point $X_k^*$ which is hyperbolic and satisfies*

$$\lim_{k\to\infty} X_k^* = X^* := \mathcal{S}(\bar{M}(Y^*)Y^*).$$

*Moreover, the following hold*:
  (i) *$X_k^*$ is asymptotically stable (resp. unstable) if and only if $Y^*$ is asymptotically stable (resp. unstable).*
  (ii) *Let $Y^*$ be asymptotically stable, and let $X_0 \in \mathbb{R}_+^{4n}$ be such that the solution $\{Y(t)\}_{t=0,1,...}$ to* (17) *corresponding to the initial data $Y_0 := UX_0$ satisfies $\lim_{t\to\infty} Y(t) = Y^*$. Then, for all $k \geq k_0$, the solution to* (15) *$\{X_k(t)\}_{t=0,1,...}$ with $X_k(0) = X_0$ satisfies $\lim_{t\to\infty} X_k(t) = X_k^*$.*

## 4 Metapopulation SEIRS epidemic model with constant recruitment, standard incidence, and constant movement rates

We consider, for each patch $j \in \{1,\ldots,n\}$, a local SEIRS model (14) assuming constant recruitment function and standard incidence for disease transmission.

$$\begin{aligned}
S_j(t+1) &= B_j + \sigma_j^S \gamma_j^R R_j(t) + \sigma_j^S \left(1 - \beta_j \frac{I_j(t)}{N_j(t)}\right) S_j(t), \\
E_j(t+1) &= \sigma_j^E \beta_j \frac{I_j(t)}{N_j(t)} S_j(t) + \sigma_j^E (1 - \gamma_j^E) E_j(t), \\
I_j(t+1) &= \sigma_j^I \gamma_j^E E_j(t) + \sigma_j^I (1 - \gamma_j^I) I_j(t), \\
R_j(t+1) &= \sigma_j^R \gamma_j^I I_j(t) + \sigma_j^R (1 - \gamma_j^R) R_j(t).
\end{aligned} \quad (18)$$

Parameters satisfy $B_j > 0$, $\beta_j \in (0,1]$ and $\sigma_j^C, \gamma_j^C \in (0,1)$, $j \in \{1,\ldots,n\}$ and $C \in \mathcal{C}$.

The associated local basic reproduction number,

$$\mathcal{R}_0^j = \frac{\sigma_j^E \sigma_j^I \gamma_j^E \beta_j}{(1 - \sigma_j^E(1-\gamma_j^E))(1 - \sigma_j^I(1-\gamma_j^I))}, \quad (19)$$

rules the disease dynamics. If $\mathcal{R}_0^j < 1$, then (Theorem 3) the DFE $X_0^j = (B_j/(1-\sigma_j^S), 0, 0, 0)$ is GAS (disease eradication), and if $\mathcal{R}_0^j > 1$, then (Theorem 4) the system is uniformly persistent (disease endemicity).

We assume that movement rates are constant. The constant regular stochastic matrices describing the movements in each compartment are $M^S, M^E, M^I, M^R \in \mathbb{R}_+^{n \times n}$. Let $M = \text{diag}(M^S, M^E, M^I, M^R) \in \mathbb{R}_+^{4n \times 4n}$ be the matrix defining the whole fast process $F(X) = MX$.

The complete two time scale model takes the form of system (15)

$$X(t+1) = \mathcal{S}(M^k X(t)), \quad (20)$$

where $\mathcal{S}$ is the map representing the disease dynamics.

### 4.1 Reduced system

We follow the procedure described in Sect. 3.1. For every $C \in \mathcal{C}$, we denote $\bar{m}^C = (m_j^C)_{j \in \{1,\ldots,n\}}$ the column right eigenvector of matrix $M^C$ associated with eigenvalue 1 that satisfies $\bar{1} \bar{m}^C = 1$, and $\bar{M} = \text{diag}(\bar{m}^S, \bar{m}^E, \bar{m}^I, \bar{m}^R)$. Thus, the reduced system associated with system (20) can be expressed in the form of (17) as

$$Y(t+1) = U\mathcal{S}(\bar{M}Y(t))$$

that corresponds to the following SEIRS model:

$$\begin{aligned}
S(t+1) &= \bar{B} + \delta_R^S R(t) + \delta_S^S S(t) - \bar{\beta}_S(Y(t)) S(t) I(t), \\
E(t+1) &= \bar{\beta}_E(Y(t)) S(t) I(t) + \delta_E^E E(t), \\
I(t+1) &= \delta_E^I E(t) + \delta_I^I I(t), \\
R(t+1) &= \delta_I^R I(t) + \delta_R^R R(t),
\end{aligned} \quad (21)$$

where the parameters are weighted means of the corresponding local parameters, with the weights being the elements of the stable probability distributions associated with the movements process:

$$\bar{B} = \sum_{j=1}^{n} B_j, \qquad \delta_R^S = \sum_{j=1}^{n} \sigma_j^S \gamma_j^R m_j^R$$

$$\delta_S^S = \sum_{j=1}^{n} \sigma_j^S m_j^S, \qquad \bar{\beta}_S(Y(t)) = \sum_{j=1}^{n} \frac{\sigma_j^S \beta_j m_j^I m_j^S}{\sum_{C \in \mathcal{C}} m_j^C C(t)},$$

$$\bar{\beta}_E(Y(t)) = \sum_{j=1}^{n} \frac{\sigma_j^E \beta_j m_j^I m_j^S}{\sum_{C \in \mathcal{C}} m_j^C C(t)}, \qquad \delta_E^E = \sum_{j=1}^{n} \sigma_j^E (1 - \gamma_j^E) m_j^E, \qquad (22)$$

$$\delta_E^I = \sum_{j=1}^{n} \sigma_j^I \gamma_j^E m_j^E, \qquad \delta_I^I = \sum_{j=1}^{n} \sigma_j^I (1 - \gamma_j^I) m_j^I,$$

$$\delta_I^R = \sum_{j=1}^{n} \sigma_j^R \gamma_j^I m_j^I, \qquad \delta_R^R = \sum_{j=1}^{n} \sigma_j^R (1 - \gamma_j^R) m_j^R.$$

System (21) is similar to the local system (18), the main difference between them being the transmission terms $\bar{\beta}_S$ and $\bar{\beta}_E$. Its recruitment constant is the sum of the local ones. Coefficients $\delta_R^S$, $\delta_S^S$, $\delta_E^E$, $\delta_E^I$, $\delta_I^I$, $\delta_I^R$, and $\delta_R^R$ are all in $(0, 1)$. The same happens to $\bar{\beta}_S(Y(t))I(t)$ and $\bar{\beta}_E(Y(t))I(t)$. In fact,

$$\bar{\beta}_S(Y(t))I(t) = \sum_{j=1}^{n} \frac{\sigma_j^S \beta_j m_j^I I(t) m_j^S}{\sum_{C \in \mathcal{C}} m_j^C C(t)} \leq \sum_{j=1}^{n} \sigma_j^S \beta_j m_j^S \leq \max_j \{\sigma_j^S \beta_j\} < 1,$$

and $\bar{\beta}_E(Y(t))I(t) \leq \sum_{j=1}^{n} \sigma_j^E \beta_j m_j^S \leq \max_j \{\sigma_j^E \beta_j\} < 1$.

We now proceed to briefly analyze system (21) following the steps carried out in Sect. 2 with system (4). As they are both very similar, we only fall into the details in case of significant difference.

All parameters are positive. Also, in the $S$ equation, $\delta_S^S - \bar{\beta}_S(Y(t))I(t)$ is positive for any $Y(t) \in \mathbf{R}_+^4$. Therefore, $\mathbf{R}_+^4$ is forward invariant under the semiflow defined by the discrete-time system (21).

**Proposition 7** *System* (21) *is dissipative.*

*Proof of Proposition 7* Let us define $\hat{\sigma} = \max_{C \in \mathcal{C}, j \in \{1,\ldots,n\}} \{\sigma_j^C\} \in (0, 1)$.

Summing up the four equations of system (21)

$$N(t+1) = \bar{B} + \big(\delta_S^S - \bar{\beta}_S(Y(t))I(t) + \bar{\beta}_E(Y(t))I(t)\big)S(t)$$
$$+ \big(\delta_E^E + \delta_E^I\big)E(t) + \big(\delta_I^I + \delta_I^R\big)I(t) + \big(\delta_R^R + \delta_R^S\big)R(t).$$

The coefficients of $S(t)$, $E(t)$, $I(t)$, and $R(t)$ are all bounded by $\hat{\sigma}$. Let us show it in the case of $S(t)$. The rest are analogous. Using the fact that the sum of a convex combination is less than or equal to the maximum of the summands, we easily obtain

$$\delta_S^S - \bar{\beta}_S(Y(t))I(t) + \bar{\beta}_E(Y(t))I(t)$$
$$= \sum_{j=1}^n \left( \sigma_j^S \left(1 - \frac{\beta_j m_j^I I(t)}{\sum_{C \in \mathcal{C}} m_j^C C(t)}\right) + \sigma_j^E \frac{\beta_j m_j^I I(t)}{\sum_{C \in \mathcal{C}} m_j^C C(t)} \right) m_j^S$$
$$\leq \sum_{j=1}^n \left( \max_j \{\sigma_j^S, \sigma_j^E\} \right) m_j^S \leq \hat{\sigma}.$$

Therefore, the total population verifies the following recurrent inequality:

$$N(t+1) \leq \bar{B} + \hat{\sigma} N(t),$$

and the proof can be completed as in Proposition 1. □

Note from the proof that, as in the case of system (4), all the nonnegative solutions of system (21) are attracted by the compact set $\{Y \in \mathbf{R}_+^4 : N \in [0, \bar{B}/(1-\hat{\sigma})]\}$.

The associated disease-free system

$$S(t+1) = \bar{B} + \delta_R^S R(t) + \delta_S^S S(t),$$
$$R(t+1) = \delta_R^R R(t)$$

is linear and, thus, it is straightforward to prove that it possesses the GAS equilibrium $(\bar{B}/(1-\delta_S^S), 0)$.

The unique DFE of system (21) is $Y_0^* = (\bar{B}/(1-\delta_S^S), 0, 0, 0)$.

The next-generation method [2] allows us to calculate its basic reproduction number

$$\overline{\mathcal{R}}_0 = \frac{\delta_E^I \bar{\beta}_I}{(1-\delta_E^E)(1-\delta_I^I)}, \tag{23}$$

where $\bar{\beta}_I = \sum_{j=1}^n \sigma_j^E \beta_j m_j^I \in (0,1)$.

The next result characterizes the eradication/endemicity of the disease in terms of $\overline{\mathcal{R}}_0$.

**Theorem 8** *Consider system* (21). *Then*
  (a) *If* $\overline{\mathcal{R}}_0 < 1$, *then DFE* $Y_0^*$ *is GAS.*
  (b) *If* $\overline{\mathcal{R}}_0 > 1$, *then DFE* $Y_0^*$ *is unstable and the system is uniformly persistent.*

*Proof of Theorem* 8
  (a) Follow the proof of Theorem 3 using that

$$E(t+1) = \delta_E^E E(t) + \left( \sum_{j=1}^n \sigma_j^E \frac{\beta_j m_j^I m_j^S S(t)}{\sum_{C \in \mathcal{C}} m_j^C C(t)} \right) I(t)$$
$$\leq \delta_E^E E(t) + \left( \sum_{j=1}^n \sigma_j^E \beta_j m_j^I \right) I(t) = \delta_E^E E(t) + \bar{\beta}_I I(t).$$

(b) As in the proof of Theorem 4, we prove, arguing by contradiction, that (21) is uniformly weakly persistent. We omit the details that are very similar in both proofs.

For any $\varepsilon > 0$, there exist $t_1 > 0$ and $\bar{r} > 0$ such that

$$E(t) + I(t) + R(t) \leq (1 + \bar{r})\varepsilon \quad \text{for } t \geq t_1. \tag{24}$$

We have the inequality $N(t+1) \geq \hat{\sigma} N(t) + \bar{B}$, where $\hat{\sigma} = \min_{C \in \mathcal{C}, j \in \{1,\ldots,n\}}\{\sigma_j^C\}$, $\hat{\sigma} \in (0,1)$, which implies that $N(t) \geq \bar{n} := \min\{N(0), B/(1-\hat{\sigma})\} > 0$.

Thus, together with (24), this yields

$$S(t) \geq \bar{n} - (1 + \bar{r})\varepsilon \quad \text{for } t \geq t_1. \tag{25}$$

A lower bound for the coefficient of $I(t)$ in the E equation, calling $\bar{m} = \max\{m_j^E, m_j^I, m_j^R\}$, is the following:

$$\bar{\beta}_E(Y(t))S(t) \geq G(\varepsilon) := \sum_{j=1}^n \frac{\sigma_j^E \beta_j m_j^I m_j^S (\bar{n} - (1+\bar{r})\varepsilon)}{m_j^S(\bar{n} - (1+\bar{r})\varepsilon) + \bar{m}(1+\bar{r})\varepsilon}.$$

Matrix

$$\bar{P}_\varepsilon = \begin{bmatrix} \delta_E^E & G(\varepsilon) \\ \delta_E^I & \delta_I^I \end{bmatrix}$$

satisfies, for $t \geq t_1$,

$$\begin{bmatrix} E(t+1) \\ I(t+1) \end{bmatrix} \geq \bar{P}_\varepsilon \begin{bmatrix} E(t) \\ I(t) \end{bmatrix}.$$

Therefore, to get the required contradiction and complete the proof, it is enough to find $\varepsilon$ such that the net reproductive rate $\overline{\mathcal{R}}_{0,\varepsilon}$ of $\bar{P}_\varepsilon$ is larger than 1. The facts that $\overline{\mathcal{R}}_{0,\varepsilon} = \frac{G(\varepsilon)}{\bar{\beta}_I} \overline{\mathcal{R}}_0$, $\overline{\mathcal{R}}_0 > 1$, and $G(\varepsilon)$ is a continuous decreasing function with $\lim_{\varepsilon \to 0^+} G(\varepsilon) = \bar{\beta}_I$ justify the existence of the required $\varepsilon$. □

A direct consequence of Theorems 6 and 8 is the following result, which we apply in the next section.

**Corollary 9** *If $\overline{\mathcal{R}}_0 < 1$, there exists an integer $k_0 \geq 0$ such that for all $k \geq k_0$ system* (20) *has an equilibrium point $X_k^*$ which is GAS and satisfies*

$$\lim_{k \to \infty} X_k^* = X^* := \mathcal{S}(\bar{M}Y_0^*).$$

We note that $X^*$ is a DFE. Indeed:

$$X^* = \mathcal{S}\left(\text{diag}(\bar{m}^S, \bar{m}^E, \bar{m}^I, \bar{m}^R)(\bar{B}/(1-\delta_S^S), 0, 0, 0)\right)$$
$$= \mathcal{S}\left(\text{diag}\left(\frac{\bar{B}}{1-\delta_S^S}\bar{m}^S, 0, 0, 0\right)\right) = \text{col}(\bar{x}_S^*, 0, 0, 0),$$

where $(\bar{x}_S^*)_j = \frac{\bar{B}m_j^S}{1-\delta_S^S}$, $j = 1, 2, \ldots, n$. Therefore, Corollary 9 states that a sufficient condition for the global disease eradication in system (20) is that the basic reproduction number of the reduced system $\overline{\mathcal{R}}_0$ is less than one.

### 4.2 Results

To illustrate the use of the developed framework, we explore the possibility of eradication of the disease through appropriate movements when it is endemic in every isolated patch.

We consider system (20). The disease is endemic in each patch $j \in \{1, \ldots, n\}$ if the associated local basic reproduction numbers satisfy

$$\mathcal{R}_0^j = \frac{\sigma_j^E \sigma_j^I \gamma_j^E \beta_j}{(1 - \sigma_j^E(1 - \gamma_j^E))(1 - \sigma_j^I(1 - \gamma_j^I))} > 1.$$

We want to find conditions such that the appropriate movements of exposed and infectious individuals drive the disease to eradication. We are looking for conditions on parameters $\sigma_j^E, \sigma_j^I, \gamma_j^E, \gamma_j^I$, and $\beta_j$, for which there exist stable equilibrium proportions of the distribution of exposed and infectious individuals, $m_j^E$ and $m_j^I$, such that

$$\overline{\mathcal{R}}_0 = \frac{\sum_{j=1}^n \sigma_j^I \gamma_j^E m_j^E \sum_{j=1}^n \sigma_j^E \beta_j m_j^I}{(1 - \sum_{j=1}^n \sigma_j^E(1 - \gamma_j^E)m_j^E)(1 - \sum_{j=1}^n \sigma_j^I(1 - \gamma_j^I)m_j^I)} < 1.$$

We make some simplifying assumptions to deal with the large number of parameters. Among the five parameters involved in the expression of $\mathcal{R}_0^j$, we consider three of them, $\sigma^E, \gamma^E$ and $\beta$, to be the same in all patches, that is,

$$\sigma_j^E = \sigma^E, \qquad \gamma_j^E = \gamma^E, \qquad \beta_j = \beta, \quad j \in \{1, \ldots, n\}.$$

Calling

$$A = \frac{\sigma^E \gamma^E \beta}{1 - \sigma^E(1 - \gamma^E)}, \qquad (26)$$

we obtain the following simplified expressions for the basic reproduction numbers:

$$\mathcal{R}_0^j = A \frac{\sigma_j^I}{1 - \sigma_j^I(1 - \gamma_j^I)} \quad \text{and} \quad \overline{\mathcal{R}}_0 = A \frac{\sum_{j=1}^n \sigma_j^I m_j^E}{1 - \sum_{j=1}^n \sigma_j^I(1 - \gamma_j^I)m_j^I}.$$

Considering $\mathcal{R}_0^j$ as a function of any of the five parameters involved in its expression, it can be checked that it is monotone on $[0, 1]$. Keeping four out of these five parameters constant across patches would make it impossible to find conditions to obtain $\overline{\mathcal{R}}_0 < 1$, since in that case

$$\overline{\mathcal{R}}_0 \geq \min_{j=1,\ldots,n} \{\mathcal{R}_0^j\} > 1.$$

Another simplification leading to the same negative conclusion is assuming coefficients $m_j^E$ and $m_j^I$ to be equal, i.e., $m_j^E = m_j^I = m_j$, since in this case $\mathcal{R}_0^j > 1$, $j \in \{1, \ldots, n\}$, implies

$A\sigma_j^I > 1 - \sigma_j^I(1 - \gamma_j^I)$, and

$$A\sum_{j=1}^n \sigma_j^I m_j > 1 - \sum_{j=1}^n \sigma_j^I(1 - \gamma_j^I)m_j,$$

which yields again $\overline{\mathcal{R}}_0 > 1$. Therefore, to obtain $\mathcal{R}_0^j > 1$ and $\overline{\mathcal{R}}_0 < 1$, we need to assume that movements lead exposed and infectious individuals to different distributions among patches.

We analyze this situation, local endemicity in isolated patches with global eradication through appropriate movements, in a simple two-patch metapopulation.

We have

$$\mathcal{R}_0^1 = A\frac{\sigma_1^I}{1 - \sigma_1^I(1 - \gamma_1^I)} \quad \text{and} \quad \mathcal{R}_0^2 = A\frac{\sigma_2^I}{1 - \sigma_2^I(1 - \gamma_2^I)}, \tag{27}$$

and, calling $x = m_1^E$ and $y = m_1^I$, the expression for $\overline{\mathcal{R}}_0$ becomes

$$\overline{\mathcal{R}}_0 = A\frac{\sigma_1^I x + \sigma_2^I(1 - x)}{(1 - \sigma_1^I(1 - \gamma_1^I))y + (1 - \sigma_2^I(1 - \gamma_2^I))(1 - y)}, \tag{28}$$

where $x$ and $y$, the fractions of exposed and infectious individuals in patch 1, take values on $(0, 1)$.

Expression $A$, depending on the values of parameters $\sigma^E, \gamma^E, \beta \in (0, 1)$, can also take any value on $(0, 1)$. On the other hand, the expression $\sigma_j^I/(1 - \sigma_j^I(1 - \gamma_j^I))$, depending on $\sigma_j^I, \gamma_j^I \in (0, 1)$, can take any positive value. Therefore, independently of the value of $A$, it is always feasible to have $\mathcal{R}_0^1 > 1$ and $\mathcal{R}_0^2 > 1$. Without loss of generality we can assume

$$1 < \mathcal{R}_0^2 \leq \mathcal{R}_0^1 \tag{29}$$

or, equivalently,

$$\frac{1}{A} < \frac{\sigma_2^I}{1 - \sigma_2^I(1 - \gamma_2^I)} \leq \frac{\sigma_1^I}{1 - \sigma_1^I(1 - \gamma_1^I)}.$$

Now, we look for the condition that ensures $\overline{\mathcal{R}}_0 < 1$, that is,

$$g(x, y) := \frac{\sigma_1^I x + \sigma_2^I(1 - x)}{(1 - \sigma_1^I(1 - \gamma_1^I))y + (1 - \sigma_2^I(1 - \gamma_2^I))(1 - y)} < \frac{1}{A},$$

where $x$ and $y$ can be chosen in $(0, 1)$. This condition can be easily expressed as

$$\min_{(x,y) \in [0,1] \times [0,1]} g(x, y) < \frac{1}{A}.$$

This minimum is attained at the boundary of $[0, 1] \times [0, 1]$ because $\partial g/\partial x$ never changes sign. In the four sides of the boundary $g$ becomes a monotone function of either variable $x$ or $y$, and therefore the minimum must be in one of the four corners. We have already assumed that $g(0, 0) = \mathcal{R}_0^2/A > 1/A$ and $g(1, 1) = \mathcal{R}_0^1/A > 1/A$. So, to get $\overline{\mathcal{R}}_0 < 1$, one of the following two inequalities must hold: either $g(1,0) < 1/A$ or $g(0,1) < 1/A$. Both cannot hold, since the first one together with $g(0,0) > 1/A$ implies that $\sigma_1^I < \sigma_1^I$, whereas the second one and $g(1,1) > 1/A$ lead to the opposite inequality.

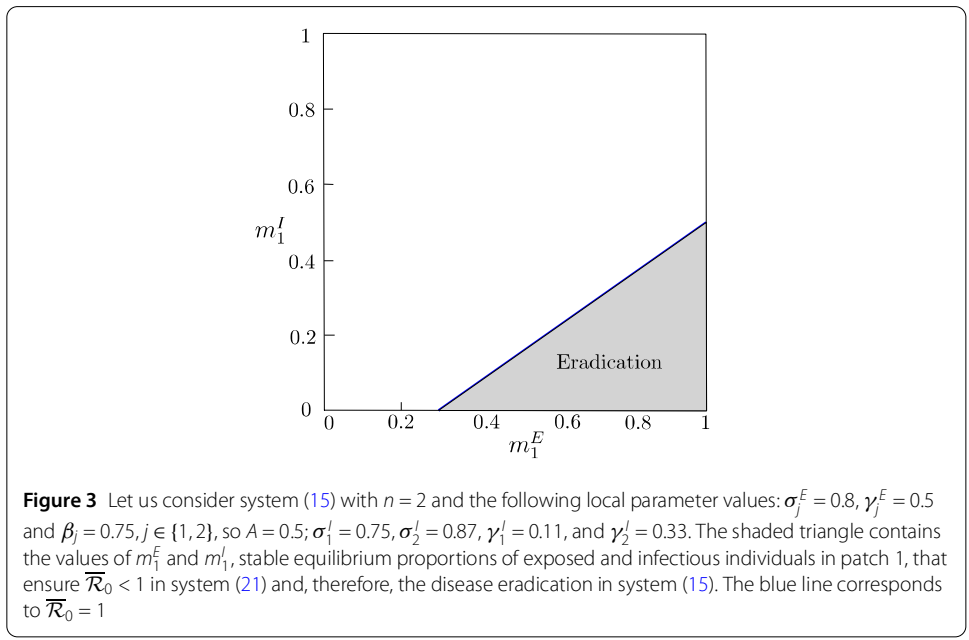

**Figure 3** Let us consider system (15) with $n=2$ and the following local parameter values: $\sigma_j^E = 0.8$, $\gamma_j^E = 0.5$ and $\beta_j = 0.75$, $j \in \{1,2\}$, so $A = 0.5$; $\sigma_1^I = 0.75$, $\sigma_2^I = 0.87$, $\gamma_1^I = 0.11$, and $\gamma_2^I = 0.33$. The shaded triangle contains the values of $m_1^E$ and $m_1^I$, stable equilibrium proportions of exposed and infectious individuals in patch 1, that ensure $\overline{\mathcal{R}}_0 < 1$ in system (21) and, therefore, the disease eradication in system (15). The blue line corresponds to $\overline{\mathcal{R}}_0 = 1$

Therefore if inequalities (29) hold and $g(1,0) < 1/A$, that is,

$$\frac{\sigma_1^I}{1 - \sigma_2^I(1 - \gamma_2^I)} < \frac{1}{A}, \qquad (30)$$

then the values $(m_1^E, m_1^I) \in (0,1) \times (0,1)$ under the line

$$\overline{\mathcal{R}}_0 = A \frac{\sigma_1^I m_1^E + \sigma_2^I(1 - m_1^E)}{(1 - \sigma_1^I(1 - \gamma_1^I))m_1^I + (1 - \sigma_2^I(1 - \gamma_2^I))(1 - m_1^I)} = 1$$

correspond to movements leading to the global eradication of the disease (see Fig. 3).

Analogously, if inequalities (29) hold and $g(0,1) < 1/A$, that is,

$$\frac{\sigma_2^I}{1 - \sigma_1^I(1 - \gamma_1^I)} < \frac{1}{A}, \qquad (31)$$

then the values $(m_1^E, m_1^I) \in (0,1) \times (0,1)$ over the line $\overline{\mathcal{R}}_0 = 1$ correspond to movements leading to the global eradication of the disease (see Fig. 4).

## 5 Discussion

Contrary to what happens in continuous time [5, 16, 17, 24, 28], in the literature there are almost no references on discrete-time epidemic models with two time scales. The authors did present some of them [9–11] to analyze the influence of a parasite on the dynamics of a community. In them, the effect of the parasite is represented by an SIS model, and it is considered fast in comparison to the community dynamics.

This work aims to establish a framework applicable to discrete-time models similar to the one existing in continuous time. However, it cannot be as flexible and powerful, for including two time scales in a discrete model has to be done sequentially and that necessarily imposes constrains. Moreover, the analysis of these models requires the reduction of their dimension, and the results at disposal to relate the behavior of the original and the reduced models are scarce [22, 23].

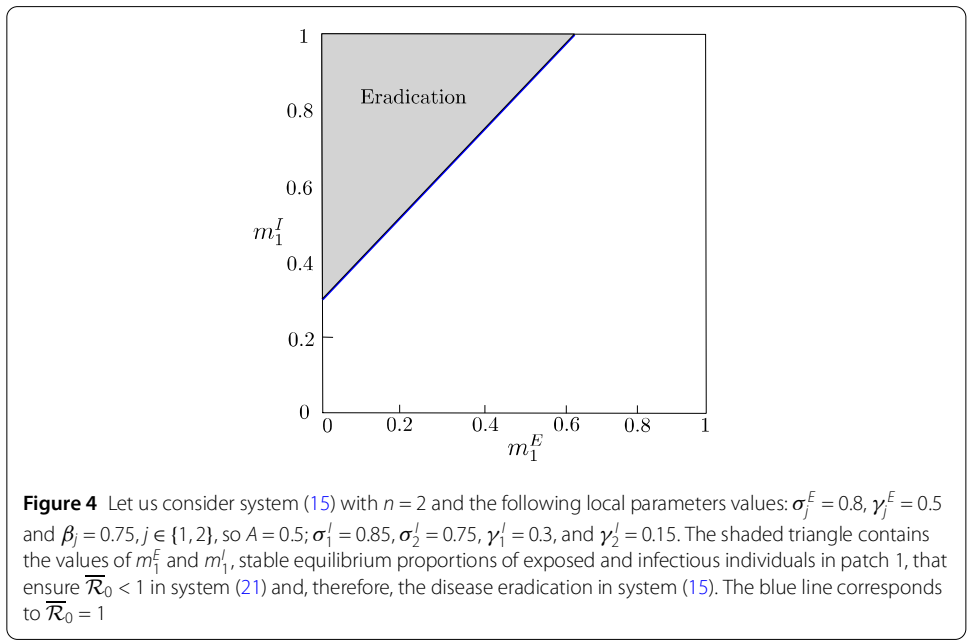

**Figure 4** Let us consider system (15) with $n = 2$ and the following local parameters values: $\sigma_j^E = 0.8$, $\gamma_j^E = 0.5$ and $\beta_j = 0.75$, $j \in \{1,2\}$, so $A = 0.5$; $\sigma_1^I = 0.85$, $\sigma_2^I = 0.75$, $\gamma_1^I = 0.3$, and $\gamma_2^I = 0.15$. The shaded triangle contains the values of $m_1^E$ and $m_1^I$, stable equilibrium proportions of exposed and infectious individuals in patch 1, that ensure $\overline{\mathcal{R}}_0 < 1$ in system (21) and, therefore, the disease eradication in system (15). The blue line corresponds to $\overline{\mathcal{R}}_0 = 1$

The proposed framework consists in a disease process acting on the slow time scale coupled with a second process whose dynamics is fast in comparison. To show the capabilities of this framework, we have tried to be fairly general in our choice of both processes. A new SEIRS discrete-time model is presented and analyzed. This model is chosen because together with its subcases (SI, SIS, SEI, SEIS, SIR, SIRS,...) it encompasses a rather large number of commonly used models. Anyway, changing the local disease model in our framework is straightforward. The second process is presented under the form of generalized individual movements between patches of a generalized metapopulation. This is the way of admitting processes as general as possible in such a way that they lead to two time scale systems susceptible of being reduced by the methods mentioned above.

The analysis of SEIRS model (4) characterizes the eradication or endemicity of the disease. Its $\mathcal{R}_0$ (10) is calculated by the next-generation method and, together with some other mild assumptions, if it is less than 1, then its DFE is GAS. On the other hand, if $\mathcal{R}_0 > 1$, the disease becomes endemic, i.e., the model is uniformly persistent in its infected compartments. When analyzing the full model with two time scales (20), the associated reduced model (21) is an SEIRS model similar to (4). Its basic reproduction number $\overline{\mathcal{R}}_0$ (23) serves to decide whether the EE is GAS or the system is persistent on infection. The main reduction result that is applied, Theorem 6, allows us to transfer the results on the asymptotic behavior of the system having to do with equilibria from the reduced system to the complete system for $k$ large enough. Thus, $\overline{\mathcal{R}}_0 < 1$ ensures the eradication of the disease in the complete model. This makes $\overline{\mathcal{R}}_0$ play the role of the basic reproduction number of model (20) with a good approximation. Therefore, the presented framework provides a method to approximate in a simpler form the basic reproduction number of a complex system.

In the general case of the complete model (15), if the reduced model (17) were to admit an endemic equilibrium (EE), possibly GAS, indicating the persistence of the disease, then Theorem 6 would apply. Thus, we could state the following corollary: If $Y_1^*$ is an EE of the reduced model (17), then, for $k$ large enough, the complete model (15) possesses an EE close to $X_1^* := \mathcal{S}(\bar{M}(Y_1^*)Y_1^*)$, and if $Y_1^*$ is GAS so is $X_1^*$. An issue that it is not covered by Theorem 6, and we will address it in a future work, has to do with the property of uniform persistence. If the reduced model (17) is uniformly persistent, Theorem 6 gives no information on the persistence of the complete model (15).

The general framework of Sect. 3 is treated in Sect. 4 for the case of constant recruitment, standard incidence, and constant movement rates. This is done to present a final pertinent illustration of the application of the process. We address the situation of a disease that is endemic in each isolated patch. Conditions on disease local parameters are found that enable driving the disease to extinction globally through a set of appropriate movements. Certainly, this application is only a minimal example of the large number of relevant aspects of disease dynamics that can be treated on similar terms using our approach.

### Acknowledgements
The authors would like to thank an anonymous referee for their suggestions to improve this work.

### Funding
Authors are supported by Ministerio de Ciencia e Innovación (Spain), Project PID2020-114814GB-I00.

### Abbreviations
DFE, Disease-free equilibrium; EE, Endemic equilibrium; GAS, Globally asymptotically stable; LAS, Locally asymptotically stable; SEIRS, Susceptible–Exposed–Infectious–Recovered-Susceptible.

### Availability of data and materials
Not applicable.

## Declarations

### Ethics approval and consent to participate
Not applicable.

### Consent for publication
Not applicable.

### Competing interests
The authors declare that they have no competing interests.

### Authors' contributions
Both authors have the same contributions. All authors read and approved the final manuscript.

### Author details
[1]U.D. Matemáticas, Universidad de Alcalá, Alcalá de Henares, Spain.  [2]Dpto. Matemática Aplicada, ETSI Industriales, Universidad Politécnica de Madrid, Madrid, Spain.

### Publisher's Note
Springer Nature remains neutral with regard to jurisdictional claims in published maps and institutional affiliations.